\documentclass[12pt, reqno]{amsart}
\usepackage{hyperref}

\usepackage{amssymb,amsmath,graphicx,amsthm}

\def\contrazione{\raisebox{1pt}{\,{\mbox{\tiny{$|\!\raisebox{-0.7pt}{\underline{\hphantom{X}}}$}}}\,}}
\def\sumKk-1{\underset{|K|=k-1}{{\sum}'}}
\def\sumKq-1{\underset{|K|=q-1}{{\sum}'}}

\def\sumJ(p-1){\underset{|J|=p-1}{{\sum}'}}

\def\sumKq{\underset{|K|=q}{{\sum}'}}
\def\sumKp-2{\underset{|K|=p-2}{{\sum}'}}

\def\sumjq+1{\underset {j\leq q+1}\sum}
\def\sumjn-1{\underset {j\leq n-1}\sum}

\def\epf{ \hskip15cm$\Box$}
\def\bpf{{\it Proof. }\hskip0.2cm}
\def\bt{\begin{theorem}}
\def\el{\end{lemma}}
\def\bl{\begin{lemma}}
\def\et{\end{theorem}}
\def\bp{\begin{proposition}}
\def\ep{\end{proposition}}
\def\bd{\begin{definition}}
\def\ed{\end{definition}}
\def\br{\begin{remark}}
\def\er{\end{remark}}

\def\simleq{\underset\sim<}

\def\simgeq{\underset\sim>}

\def\T{\text}

\def\1#1{\overline{#1}}

\def\2#1{\widetilde{#1}}

\def\3#1{\widehat{#1}}

\def\4#1{\mathbb{#1}}

\def\5#1{\frak{#1}}

\def\6#1{{\mathcal{#1}}}

\def\C{{\4C}}

\def\R{{\4R}}

\def\sumJ{\underset{|J|=k}{{\sum}'}}

\def\sumjq{\underset {j\leq q}\sum}

\def\T{\text}
\newcommand{\Om}{\Omega}
\newcommand{\om}{\omega}

\newcommand{\no}[1]{\|{#1}\|}

\newcommand{\No}[1]{\|{#1}\|^2}
\def\NO#1{||#1||^2}
\def\R{{\Bbb R}}

\def\C{{\Bbb C}}

\def\di{\partial}
\def\dib{\bar\partial}
\def\Label#1{\label{#1}}

\def\Op{\T{Op}^{\T{ord}(\Psi)-\frac12}}
\def\Opbis{\T{Op}^{2\,\T{ord}(\Psi)-1}}

\def\simleq{\underset\sim<}

\def\simgeq{\underset\sim>}

\def\T{\text}

\def\1#1{\overline{#1}}

\def\2#1{\widetilde{#1}}

\def\3#1{\widehat{#1}}

\def\4#1{\mathbb{#1}}

\def\5#1{\frak{#1}}

\def\6#1{{\mathcal{#1}}}

\def\C{{\4C}}

\def\R{{\4R}}

\def\sumJ{\underset{|J|=k}{{\sum}'}}

\def\sumjq{\underset{j=1}{\overset{q_0}\sum}}

\numberwithin{equation}{section}

\def\T{\text}

\theoremstyle{plain}

\newtheorem{theorem}{Theorem}[section]

\newtheorem{lemma}[theorem]{Lemma}

\newtheorem{proposition}[theorem]{Proposition}

\theoremstyle{definition}

\newtheorem{definition}[theorem]{Definition}

\theoremstyle{remark}

\newtheorem{remark}[theorem]{Remark}

\begin{document}

\title[The Kohn-H\"ormander-Morrey...]{ The  Kohn-H\"ormander-Morrey formula twisted by a pseudodifferential operator}         
\author[L.Baracco, M.Fassina and S.Pinton]{Luca Baracco, Martino Fassina and Stefano Pinton}        
\address{Dipartimento di Matematica, Universit\`a di Padova, via 
Trieste 63, 35121 Padova, Italy}
\email{baracco@math.unipd.it,fassina@math.unipd.it,pinton@math.unipd.it}

\maketitle
\begin{abstract}
We establish a general, weighted Kohn-H\"ormander-Morrey formula twisted by a pseudodifferential operator. As an application, we exhibit a new class of domains for which the $\dib$-Neumann problem is locally hypoelliptic.
\end{abstract}
\tableofcontents 
\section{Introduction}
\Label{s1}
The Kohn-H\"ormander-Morrey formula, together with the density of $C^\infty$- into $L^2$-forms in the graph norm yields the closure of the  range of $\dib$ and $\Box$  \cite{KN65}, \cite{FK72}.
By adding a weight, it also gives the global boundary regularity of the (non-canonical) solution of $\dib$. When it comes with a ``gain" it gives in fact the  regularity of the canonical solution of the $\dib$-Neumann problem \cite{K73}. For a gain consisting of the multiplication by a big constant, as in the compactness estimates, the regularity is global. For a gain such a subelliptic \cite{KN65} or superlogarithmic \cite{K02} the regularity is stronger, that is, local: $u$ is regular precisely in the portion of $b\Om$  where $\Box u $ is. When these latter  gains do not occur, but the points of failure are confined to a real curve transversal to the CR directions, local regularity still holds \cite{K00} and \cite{BKZ14}. This is an exquisitely geometric conclusion. 
In \cite{BPZ14} it is  shown that   good  estimates in full are not needed and what really counts is that for a system of cut-off $\{\eta\}$,  the gradient $\di\eta$ and the Levi form $\di\dib\eta$
are good multipliers in the sense of Kohn \cite{K79}. If these are subelliptic multipliers, then $\Box$  is hypoelliptic. The proof consists in modifying the Kohn-H\"ormander-Morrey formula by a weight $\phi=t|z|^2-\log\eta^2$; the exploitation  of $t|z|^2$ is usual in controlling the commutators $[\dib,\Lambda^s]$ and $[\dib^*,\Lambda^s]$, but the one of $-\log\eta^2$ is new and is designed to avoid the commutators $[\dib,\eta]$ and $[\dib^*,\eta]$. 
Inserting the cut-off in the weight is not the only way to proceed and this could appear as well as a twisting term like in the formula of \cite{S10}; what is crucial is not to apply the basic estimate to $\eta u$, for a testing form $u$, but the formula twisted by $\eta$ to the plain $u$. 
If the multipliers are weaker than subelliptic one needs a stronger modification of the Kohn-H\"ormander-Morrey formula in which not only the cut-off but also a general pseudodifferential operator appears as already commutated with $\dib$ and $\dib^*$. 
The motivation of the present paper is to present, for the $\dib$-Neumann problem, a general  Kohn-H\"ormander-Morrey formula with weight in which    a general pseudodifferential operator appears as a twisting term. This has already been done in \cite{BKPZ14} for the tangential system.
The present paper serves therefore as the $\dib$-Neumann version of the tangential twisted estimate established in \cite{BKPZ14}.
It proves useful for the local regularity of $\Box$  when the multipliers are not subelliptic but only superlogarithmic (in the slightly stronger sense of \eqref{multiplier} below). This requires to twist the formula not only  by $\eta$ but also by $R^s$ where $R^s$ is the modification of the standard elliptic operator  of order $s$ introduced by Kohn in \cite{K02}. As an application we get a new class of domains of infinite type for which $\Box$ is locally hypoelliptic.

\section{The twisted basic estimate for the $\dib$-Neumann problem}
\Label{s2}
Let $\Om$ be a domain of $\C^n$ with $C^\infty$-boundary $b\Om$, $z_o$ a point of $b\Om$, $U$ an open neighborhood of $z_o$ and $\Psi$ a tangential pseudodifferential operator with real symbol $\mathcal S(\Psi)$. This is defined by introducing a local  straightening  $b\Om\simeq \R^{2n-1}\times\{0\}$ and $\Om\simeq \R^{2n-1}\times\R^-_r$ for a defining function $r<0$ of $\Om$ with $\di r\neq0$, taking coordinates $(x,r)$ or $(y,r)$ in $U$,  dual coordinates $\xi$ of $x$ and setting
$$
\Psi(u)=\int e^{ix\xi}\mathcal S(\Psi)(x,\xi)\left(\int e^{-iy\xi} u(y)dy\right)d\xi.
$$
One of the most common choice of $\Psi$ is the elliptic standard operator of degree $s$ with symbol
$$
\mathcal S (\Lambda^s)=\Big(1+|\xi|^2\Big)^{\frac12}.
$$
It is also of great interest its local modification by means of a cut-off $\sigma\in C^\infty_c(U)$, which is $1$ in a neighborhood of $z_o$,
\begin{equation}
\Label{nova}
\mathcal S(R^s)=\Big(1+|\xi|^2\Big)^{s\sigma(x)}.
\end{equation}
The operators we have in mind are $\Psi=\eta\Lambda^s\eta_o$ for a pair of cut-off $\eta_o\prec \eta$ in $C^\infty_c(U)$ or $\Psi=\eta R^s\eta$ for $\eta_o\prec\sigma\prec \eta$. However, our formula applies to a general symbol and what we ask is that this is already subject to the multiplication by a cut-off, that is, $\Psi=\eta\Psi$ for some $\eta\in C^\infty_c(U)$.

To give our formula a bigger flexibility, we also consider  spaces $L^2_\phi$ weighted by a weight $e^{-\phi}$ with norm defined by $\NO{u}_\phi=\int e^{-\phi}|u|^2 dV$. In particular, the action of $\Psi^2=\eta^2$ can be achieved by means of the weight $\phi=-\log\eta^2$ and it is under this appearence that the basic formula was introduced in \cite{BPZ14}. We use the notation
$$
Q^\phi_\Psi(u,u)=\NO{\Psi \dib u}_\phi+\NO{\Psi \dib^* u}_\phi,\quad u  \in   D_{\dib}\cap {D_{\dib^*}}.
$$
We choose a smooth orthonormal basis of $(1,0)$ forms $\om_1,...,\om_{n-1},\,\om_n=\di r$ and the dual basis of vector fields $\di_{\om_1},...,\di_{\om_n}$; to simplify notation, we also write $\di_j$ instead of $\di_{\om_j}$. We use the notation $\phi_i$ for $\di_i\phi$, $\phi_{\bar i}$ for $\bar\di_i\phi$ and $(\phi_{ij})$ for the matrix of the Levi form $\di\dib\phi$. Note that $\phi_{ij}$ differs from $\di_i\dib_j(\phi)$ because of the presence of the derivatives of the coefficients of the forms $\bar\di_j$.
We  define various constants $c_{ij}=c_{ij}^n$ and $c_{ij}^h,\,\,i,j,h=1,...,n-1$ by means of the identity
\begin{equation}
\Label{2.5}
[\di_i,\dib_j]=c_{ij}(\di_n-\dib_n)+\sum_{h=1}^{n-1}c_{ij}^h\di_h-\sum_{h=1}^{n-1}\bar c_{ji}^h\dib_h.
\end{equation}
Thus $(c_{ij})$ is the matrix of the Levi form  $\di\dib r|_{T^\C b\Om}$ in the basis $\{\om_j\}$. (Here we have used the notation  $T^\C b\Om=Tb\Om\cap iTb\Om$.)
We denote by $\Op$, resp. $\T{Op}^0$, an operator of order $\T{ord}(\Psi)-\frac12$, resp. $0$. We assume that $\Op$ has symbol with support in a neighborhood of that of $\Psi$ and that $\T{Op}^0$ only depends on the $C^2$-norm of $b\Om$ and not on $\phi$ nor $\Psi$. 
Here is the substance of the paper
\bt
\Label{t2.1}
For any $u\in D^k_{\dib}\cap D^k_{\dib^*}\cap C^\infty_c(U\cap \bar\Omega),\,k\in[1,n-1]$, we have
\begin{equation}
\Label{2.1}
\begin{split}
\int_{b\Om}e^{-\phi}(c_{ij})&(\Psi u,\overline{\Psi u})dV+\int_\Om e^{-\phi}\di\dib\phi(\Psi u,\overline{\Psi u})dV+\sum_{j-1,...,n}\NO{\bar L_ju}_\phi\simleq Q^\phi_{\Psi}(u,\bar u)\\&+\left|\int_\Om e^{-\phi}[\di,[\dib,\Psi^2]](u,\bar u)dV\right|
+\NO{[\di,\phi]\contrazione \Psi u}_\phi+\NO{[\di,\Psi]\contrazione u}_\phi
\\
&+\Big|\sum_h\int_\Om (c_{ij}^h)\Big([\di_{h},\Psi](u),\bar u\Big)\,dV\Big|+Q^\phi_{\Op}(u,\bar u)+\NO{\Op u}_\phi+\NO{\Psi u}_\phi.
\end{split}
\end{equation}
\et
\br
In our application in Section~\ref{s3} below, $[\di,[\dib,\Psi^2]]$ and $[\di,\Psi]\contrazione$ have good estimates. Also, $\phi$ has ``selfbounded gradient", that is
$$
[\di,\phi]\contrazione<\di\dib \phi,
$$
where inequality is meant in the operator sense. In particular, the  term in the right of \eqref{2.1} which involves $[\di,\phi]\contrazione$  is absorbed in the left.
\er
\br
Formula \eqref{2.1} is also true for complex $\Psi$. In this case, one replaces $[\di,[\dib,\Psi^2]]$ by $[\di,[\dib,|\Psi|^2]]$ and add an additional error term $[\di,\bar\Psi]\contrazione$.
\er
\bpf
We start from 
\begin{equation}
\Label{2.2}
e^{\phi}\Psi^{-2}[\dib_i,e^{-\phi}\Psi^2]=-\phi_{\bar i}+2\frac{[\dib_i,\Psi]}\Psi+\frac{\Opbis}{\Psi^2},
\end{equation}
whose sense is fully clear when both sides are multiplied by $\Psi^2$.
In other terms, we have
\begin{equation}
\Label{2.4}
\dib^*_{e^{-\phi}\Psi^2}=\dib^*+\di\phi\contrazione-2\frac{[\di,\Psi]}\Psi\contrazione+\frac{\Opbis}{\Psi^2}.
\end{equation}
We are thus lead to define the transposed operator to $\dib_i$ by
\begin{equation}
\Label{2.3}
\delta_i:=\di_i-\phi_i+2\frac{[\di_i,\Psi]}\Psi+\frac{\Opbis}{\Psi^2}+\T{Op}^0.
\end{equation}
 Using the trivial identity $\di\dib=-\dib\di$, we have
\begin{equation}
\Label{2.6}
\begin{split}
[\delta_i&,\dib_j]=[\di_i,\dib_j]+\phi_{ij}-\sum_{h=1}^{n}c_{ij}^h\phi_h-2\frac{[\di_i,[\dib_j,\Psi]]}\Psi+2\frac{[\di_i,\Psi]\otimes[\dib_j,\Psi]}{\Psi^2}+\frac{\Opbis}{\Psi^2}+\T{Op}^0
\\
&=c_{ij}(\delta_n-\dib_n)-\sum_{h=1}^{n-1}\bar c_{ji}^h\dib_h+\sum_{h=1}^{n-1} c_{ij}^h\delta_h+\phi_{ij}-2\sum_hc^h_{ij}\frac{[\di_h,\Psi]}{\Psi}\contrazione
\\
&-2\frac{[\di_i,[\dib_j,\Psi]]}\Psi+2\frac{[\di_i,\Psi]\otimes[\dib_j,\Psi]}{\Psi^2}+\frac\Opbis{\Psi^2}+\T{Op}^0.
\end{split}
\end{equation}
We also have to observe that (cf. \cite{BKPZ14})
$$
\NO{\Psi\dib^{(*)}u}_\phi=\int_\Om e^{-\phi}\Psi^2|\dib^{(*)}u|^2\,dV+Q_{\Op}(u,\bar u)+\NO{\Op u}_\phi+\NO{\T{Op}^0\Psi u}_\phi.
$$
This yields the ``basic estimate with weight $e^{-\phi}\Psi^2$"
\begin{multline}
\Label{2.7}
\int_{b\Om}e^{-\phi}(c_{ij})(\Psi u,\overline{\Psi u})\,dV+\int_\Om[\di,[\dib,e^{-\phi}\Psi^2]](u,\bar u)\,dV
\\
-\NO{[\di,\phi]\contrazione \Psi u}_\phi-\NO{[\di,\Psi]\contrazione  u}_\phi
+\sum_j\NO{\Psi\bar\di_j u}_\phi
\\
\simleq \NO{\Psi\dib u}_\phi+\NO{\Psi\dib^*_{e^{-\phi}\Psi^2}u}_\phi+sc\NO{\Psi\bar\nabla u}_\phi+\Big|\sum_h\int_\Om (c_{ij}^h)\Big([\di_{h},\Psi](u),\bar u\Big)\,dV\Big|
\\+Q^\phi_{\Op}(u,\bar u)+\NO{\Op u}_\phi+\NO{\Psi u}_\phi.
\end{multline}
In \eqref{2.7} we absorb the term which comes with sc and rewrite $[\di,[\dib,e^{-\phi}\Psi^2]]$ by the aid of \eqref{2.2}; what we get is
\begin{multline}
\Label{2.8}
\int_{b\Om}e^{-\phi}\Psi^2 c_{ij}(u,\bar u)dV+\int_{\Om}e^{-\phi}\Psi^2 \phi_{ij}(u,\bar u)dV-\NO{[\di,\phi]\contrazione \Psi u}_\phi
\\+\int_\Om e^{-\phi}\Psi^2[\di_i,[\dib_j,\Psi]](u,\bar u)dV
-\NO{[\di,\Psi]\contrazione u}_\phi+\NO{\Psi\bar\nabla  u}_\phi
\\
\simleq \NO{\Psi\dib u}_\phi+\NO{\Psi\dib^*_{e^{-\phi}\Psi^2}u}_\phi
+Q^\phi_{\Op}(u,\bar u)+\Big|\sum_h\int_\Om (c_{ij}^h)\Big([\di_{h},\Psi](u),\bar u\Big)\,dV\Big|
\\
+\NO{\Op u}_\phi+\NO{\Psi u}_\phi.
\end{multline}
To carry out our proof we need to replace $\dib^*_{e^{-\phi}\Psi^2}$ by $\dib^*$. We have from \eqref{2.4}
\begin{equation}
\Label{2.9}
\begin{split}
\NO{\Psi \dib^*_{e^{-\phi}\Psi^2}u}_\phi&\simleq \NO{\Psi\dib^* u}_\phi+\NO{\Psi\di \phi\contrazione u}_\phi+\NO{[\dib,\Psi]\contrazione u}_\phi
\\
&+\underset{\T{\#}}{\underbrace{2\Big|\Re e (\Psi\dib^* u,\overline{\Psi\di\phi\contrazione u})_\phi\Big|+2\Big|\Re e (\Psi\dib^* u,\overline{[\di,\Psi]\contrazione u})_\phi+2\Big|\Re e (\Psi\di\phi\contrazione u,\overline{[\di,\Psi]\contrazione u})_\phi\Big|}}.
\end{split}
\end{equation}
We next estimate by Cauchy-Shwarz inequality
\begin{equation}
\Label{supernova}
\#\simleq \NO{\Psi\dib^* u}_\phi+\NO{\Psi\di\phi\contrazione u}_\phi+\NO{[\di,\Psi]\contrazione u}_\phi.
\end{equation}
We move the third, forth and fifth terms from the left to the right of \eqref{2.8}, use \eqref{2.9} and \eqref{supernova} and end up with \eqref{2.1}. 

\epf

\section{$F$ type, twisted $f$ estimate and hypoellipticity of $\Box$.}
\Label{s3}
We start by recalling a result by \cite{KZ10}. In our presentation it contains a specification of the estimate by the Levi form which is important for further  application. We consider a  smoothly bounded pseudoconvex domain $\Om\subset\C^n$. For a form $u\in D_{\dib^*}$, let $u=\sum_k\Gamma_ku$ be the decomposition into wavelets (cf. \cite{K02}), and $Q(u,\bar u)=\NO{\dib u}+\NO{\dib^*u}$   the energy. We use the notation $(c_{ij})$  and $(\phi_{ij})$ for the Levi form of the boundary 
$b\Om$ and of a function $\phi$ respectively. We introduce a real function $F$ such that $\frac{F(d)}{d}\searrow0$ as $d\searrow0$ and set $f(t):=(F^*(t^{-\frac12}))^{-1}$ where $F^*$ denotes the inverse.
\bt 
\Label{t3.1}
Assume that $b\Om$ has type $F^2$ along a submanifold $S\subset b\Om$ of CR dimension $0$ in the sense that $(c_{ij})\simgeq  \Big(\frac{F(d_S)}{d_S}\Big)^2\T{Id}$ where $d_S$ is the Euclidean distance to $S$ and Id the identity of $T^\C b\Om$. Then there is a uniformly bounded family of weights $\{\phi^k\}$ which yield the ``$f$ estimate"
\begin{equation}
\Label{3.1}
\begin{split}
\NO{f(\Lambda)u}_0&\simleq \int_{b\Om}(c_{ij})(u,\bar u)\,dV+\sum_k\int_\Om (\phi^k_{ij})(\Gamma_ku,\overline{\Gamma_k u})+\sum_{j-1,...,n}\NO{\bar L_ju}_\phi+\NO{u}_0
\\
&\simleq Q^\Om(u,\bar u)\quad\T{for any $u\in D^k_{\dib}\cap D^k_{\dib^*}\cap C^\infty_c(U\cap \bar\Omega),\,k\in[1,n-1]$.}
\end{split}\end{equation}
\et
The estimate $\NO{f(\Lambda)u}_\Om\simleq Q_\Om(u,\bar u)$ is stated in \cite{KZ10} as a combination of Theorems 1.4 and 2.1 therein. For the purpose of the present paper (Theorem \ref{t3.3} below), the two separate estimates which occur in \eqref{3.1}, with the intermediate term which carries the Levi forms, are essential.

\bpf
We give two parallel proofs inspired to \cite{KZ10}, resp. \cite{BKPZ14}, which use the families of weights
$$
\psi^k:=-\log(\frac{-r}{2^{-k}}+1)+\chi(\frac{d_S}{a_k})\log(\frac{d_S^2}{a_k^2}+1)\quad\T{resp. }\phi^k:=\chi(\frac{d_S}{a_k})\log(\frac{d_S^2}{a_k^2}+1).
$$
Here $r=0$ is an equation for $b\Om$ with $r<0$ on $\Om$, $a_k:=F^*(2^{-\frac k2})$ and $\chi$ is a cut-off such that $\chi\equiv1$ in $[0,1]$ and $\chi\equiv0$ for $s\ge2$. 
We also use the notation $S_{a_k}$ for the strip $S_{a_k}:=\{z\in\Om:\,d_{b\Om}(z)<a_k\}$. 
Following word by word the proof of \cite{KZ10}, resp. \cite{BKPZ14}, we conclude
$$
\NO{f(\Lambda)u}_0\simleq\sum_k\int_{S_{2a_k}\setminus S_{a_k}}(c_{ij})(d_{b\Om}^{-\frac12}u,\overline{d_{b\Om}^{-\frac12}u})\,dV+\sum_k\int_\Om (\phi^k_{ij})(\Gamma_k u,\overline{\Gamma_ku})+\NO{u}_0
$$
resp.
$$
\NO{f(\Lambda)u}_0\simleq\int_{\Om}(c_{ij})(\Lambda^{\frac12}u,\overline{\Lambda^{\frac12}u})\,dV+\sum_k\int_\Om (\phi^k_{ij})(\Gamma_k u,\overline{\Gamma_ku})+\NO{u}_0.
$$
Finally, from
$$
\no{d_{b\Om}^{-\frac12}u}_0\simleq \no{u}^b_0+\sum_{j=1}^n\no{\bar L_j u}_0,\quad\T{ resp. }\no{\Lambda^{\frac12}u}_0\simleq \no{u}^b_0+\sum_{j=1}^n\no{\bar L_j u}_0,
$$
(cf. \cite{K02} Section~8) we get the first estimate in \eqref{3.1}. The second is a basic estimate weighted by the $\phi^k$'s in which the weights has been removed from the norms on account of their uniform boundedness.

\epf

We modify the weights $\phi^k$ to $\phi^k+t|z|^2$ so that their Levi form releases an additional  $t\T{Id}$ for $t$ big. They are absolutely uniformly bounded with respect to $k$ and to $t$ provided that we correspondingly shrink the neighborhood $U=U_t$.  Possibly by raising to exponential, boundedness implies ``selfboundedness of the gradient" when the weight is plurisubharmonic. In our case, in which to be positive is not $(\phi^k_{ij})$ itself but $2^k(c_{ij})+(\phi^k_{ij})$, we have, for $|z|$ small
\begin{equation}
\Label{stranova}
\begin{split}
|\di_b\phi|^2&=|\di_b(\phi^k+t|z|^2)|^2
\\
&\simleq |\di_b\phi^k|^2+t^2|z|^2
\\
&\leq 2^k(c_{ij})+(\phi^k_{ij})+t.
\end{split}
\end{equation}
Going back to \eqref{2.1} under this choice of $\phi$, we have that $\no{\di_b\phi\contrazione  \Psi u}^2$ can be removed from the right side. We combine Theorem~\ref{2.1} with Theorem~\ref{t3.1} formula \eqref{3.1}, observe again that the weights $\phi^k$ can be removed from the norms by uniform boundedness,  and get the proof of the following

\bt
\Label{t3.2}
Let $b\Om$ have type $F^2$ along $S$ of CR dimension 0. Then we have the $f$ estimate
\begin{equation}
\Label{3.2}
\begin{split}
||f(\Lambda)\Psi u||^2_0&\simleq \int_{b\Om} (c_{ij})(\Psi u,\overline{\Psi u})\,dV+\sum_k\int_\Om(\phi^k_{ij})(\Gamma_k\Psi u,\overline{\Gamma_k\Psi u})\,dV+\sum_j\NO{\bar L_j\Psi u}_0+t\NO{\Psi u}_0
\\
&\simleq Q_\Psi(u,\bar u)+\NO{[\di,\Psi]\contrazione u}_0+\Big|\int_\Om [\di,[\dib,\Psi^2]](u,\bar u)\,dV\Big|
\\
&+\Big|\sum_h\int_\Om (c_{ij}^h)\Big([\di_{h},\Psi](u),\bar u\Big)\,dV\Big|+Q^\phi_{\Op}(u,\bar u)+\NO{\Op u}_0+\NO{\Psi u}_0.
\end{split}
\end{equation}
\et
We have as application a criterion of regularity for the Neumann operator $N$ in a new class of domains.  Let $b\Om$ be a ``block decomposed", rigid,  boundary, that is, defined by $x_n=\sum_{j=1}^m h^{I^j}(z_{I^j})$ where $z=(z_{I^1},...,z_{I^m},z_n)$ is a decomposition of  coordinates.
\bt
\Label{t3.3}
Assume that
\begin{equation}
\Label{3.4}
\begin{cases}
 \T{ (a) $h^{I^j}$ has infraexponential type along a totally real $S^{I^j}\setminus\{0\}$}
\\\hskip0.5cm\T{ where $S^{I^j}$ is  totally real  in $\C^{I^j}$,}
\\
 \T{ (b) $h^{I^j}_{z_j}$ are superlogarithmic multipliers at $z_{I^j}=0$.}
\end{cases}
\end{equation}
Then,  we have local hypoellipticity of $\Box$  at $z_o=0$. 
\et
In the same class of domains, it is proved in \cite{BKPZ14} the hypoellipticity of the Kohn-Laplacian $\Box_b$. 

\bpf
The proof is the same as in \cite{BKPZ14} Theorem~1.11. The argument is that, for a system  $\{\eta\}$ of cut-off, and with the decomposion $\di\eta=(\di_\tau\eta,\di_\nu\eta)$, the vectors $\di_\tau\eta$  are superlogarithmic multipliers in the sense that
for any $\epsilon$, suitable $c_\epsilon$ and for a bounded family of weights $\phi^k$, we have
\begin{equation}
\Label{multiplier}
\NO{\log(\Lambda)\di_\tau\eta\contrazione u_\tau}\le \epsilon\Big(\int_{b\Om} (c_{ij}(u_\tau,\overline{ u_\tau})dV+\sum_k\int_\Om(\phi_{ij}^k)(\Gamma_ku_\tau,\overline{\Gamma_ku_\tau})dV\Big)+c_\epsilon\NO{u_\tau}_0.
\end{equation}
This is an immediate consequence of the hypotheses (a) and (b) of Theorem~\ref{t3.3} combined with \eqref{3.2} for $f=\epsilon^{-1}\log$ (with an error $c_\epsilon$) and $\Psi=\T{id}$.  (Note that this notion of superlogarithmic multiplier is a little more restrictive than in the literature where the term between brackets in the right of  \eqref{multiplier} is replaced by $Q$.)
 $\di_\nu\eta$ is $1$ but it hits $u_\nu$ which is 0 at $b\Om$ and therefore enjoys elliptic estimates.
In the same way, $[\di_\tau,[\di_{\bar\tau},\eta]]$ are superlogarithmic multipliers, whereas $[\di_\nu,[\di_{\bar\tau},\eta]]$ and $[\di_\nu,[\di_{\bar\nu},\eta]]$ give $\frac12$-subelliptic and elliptic estimates respectively. 
We also take $\eta_o\prec\sigma\prec\eta$ and recall the operator $R^s$ whose symbol has been defined by \eqref{nova}; note that $\di_\tau\sigma$ is a superlogarithmic multiplier.
We also have to notice that,  $b\Om$ being rigid, then $(c_{ij}^h)\simleq (c_{ij})$ are  $\frac12$ subelliptic matrix multipliers (cf. \eqref{2.5} which defines these matrices).
We also have to remark that
\begin{equation}
\Label{*}
[\di,\eta R^s \eta]=\eta \di\sigma\log(\Lambda)R^s\eta+\T{Op}^{-\infty}.
\end{equation}
After this preliminary, 
and under the choice $\Psi=\eta R^s\eta$, we have readily the following chain of estimates
\begin{equation}
\begin{split}
t\NO{\Lambda^s\eta_ou}_0&\simleq t\NO{\eta R^s\eta u}_0+\NO{u}_0
\\
&\underset{\T{since $(c_{ij})\ge0$}}\simleq \Big(\int_{b\Om}(c_{ij})(\eta R^s\eta u,\overline{\eta R^s\eta u})\,dV+\sum_{k=1}^{+\infty}\int_{\Om}(\phi_{ij}^k)(\Gamma_k\eta R^s\eta u,\overline{\Gamma_k\eta R^s\eta u})\,dV\Big)+t\NO{\eta R^s\eta u}_0
\\
&\underset{\T{\eqref{2.1}}}\simleq Q^\Om_{\eta R^s\eta}(u,\bar u)+\NO{[\di ,\eta R^s\eta]\contrazione u}_0+\Big|\int_{\Om}[\di ,[\dib ,\eta R^s\eta]](u,\bar u)\,dV\Big|
\\&\hskip0.3cm+\Big|\sum_h\int_\Om (c_{ij}^h)([\di_{\om_h},\eta R^s\eta](u),\overline{\eta R^s\eta  u})\,dV\Big|+Q^\Om_{\Op}(u,\bar u)+\NO{\Op u}_0
\\
&
\underset{\overset{\T{by \eqref{*} and because $(c^h_{ij})$}}{\T {are subelliptic multipliers}}}\simleq  Q^\Om_{\eta R^s\eta}(u,\bar u)+\NO{\di (\sigma)\contrazione\log(\Lambda)\eta R^s\eta u}_\Om+Q_{\Lambda^{s-\frac12}\eta'}^\Om(u,\bar u)+\NO{\eta'u}_{s-\frac12}
\\
&\underset{\overset{\T{\eqref{multiplier}and elliptic}}{\T{ estimate for $u^\nu$}}}\simleq Q^\Om_{\eta R^s\eta}(u,\bar u)+\epsilon\Big(\int_{b\Om}(c_{ij})(\eta R^s\eta u,\overline{\eta R^s\eta u})\,dV+\sum_k\int \Big((\phi^k_{ij})(\eta R^s\eta \Gamma_ku,\times
\\
&\times \overline{\eta R^s\eta \Gamma_ku})\Big)\,dV\Big)
+c_\epsilon\NO{\eta R^s\eta u}_0+Q_{\Lambda^{s-\frac12}\eta'}^\Om(u,\bar u)+\NO{\eta'u}_{s-\frac12}
\\
&\underset{\overset{\T{absorbtion in the second line}}{\T{for $\epsilon\simleq \frac1t$ and $t\sim c_\epsilon^{-1}$}}}\simleq Q^\Om_{\eta R^s\eta}(u,\bar u)+Q_{{\Lambda^{s-\frac12}\eta'}}^\Om(u,\bar u)+(\NO{\eta'u})^\Om_{s-\frac12}.
\end{split}
\end{equation}
Now, the $s-\frac12$ norm is reduced to $0$ by induction and the various $Q^\Om_{\eta R^s\eta}$ and $Q^\Om_{{\Lambda^{s-\frac12}\eta'}}$ are estimated by a common $Q^\Om_{\Lambda^s\eta}$ for a new $\eta$. Thus we end up with
\begin{equation}
\Label{local}
\no{\eta_ou}_s\simleq \no{\eta\dib u}_s+\no{\eta\dib^*u}_s+\no{u}_0.
\end{equation}
We observe that \eqref{local} is an a-priori estimate, that is, it only holds in principle for smooth forms. To overcome this restraint,
we use the Kohn microlocal decomposition $u=u^++u^-+u^0$; since the $\dib$-Neumann problem is elliptic for $u^-$ and $u^0$, we only have to prove Theorem~\ref{t3.3} for $u^+$.   We then use an approximation of the identity in the $y_n$-variable  by smooth functions $\chi_\nu(y_n)$ and smoothen $u^+$ by $u^+_\nu:=u^+*\chi_\nu$. Now, since $b\Om$ is rigid, then $\dib^{(*)}(u^+*\chi_\nu)=(\dib^{(*)}u)^+*\chi_\nu+\tilde u^0$, where $\tilde u^0$ is an error term supported by the elliptic microlocal region. Then, from \eqref{local} we get the following. If $\dib u,\,\,\dib^*u$ belong to $H^s$ in a neighborhood of supp$\,\eta$, we have that $\eta_o u\in H^s$ (and, moreover, $\No{\eta_o u}_s\simleq Q_{\eta R^s\eta}(u,u)+\NO{u}_0$). In particular, $u\in H^s$ in a neighborhood of $\{z:\,\eta_o(z)\equiv1\}$. This concludes the proof of the theorem. 

\epf

In case of a single block $x_n=h^{I^1}$ we regain \cite{BPZ14} which   transfers \cite{K00} from the tangential system to the $\dib$-Neumann problem and also gives a more general statement. The proof is far more efficient because it uses the elementary decomposition $\di\eta=(\di_b\eta,\di_\nu\eta)$ instead of $Q=Q^\tau\oplus \bar L_n$ (over tangential forms $u^\tau$) which requires the heavy technicalities of the harmonic extension.  

\br
\Label{r3.1}
the above proof shows a general criterion. If $\di_b\eta$ is a superlogarithmic multiplier for $\Box_b$, then it is also for the $\dib$-Neumann problem; thus we have hypoellipticity of $\Box_b$ and $\Box$ at the same title.
\er

\noindent
{\it Example} 
Let $b\Om$ be defined by
$$
x_n=\sum_{j=1}^{n-1}e^{-\frac1{|z_j|^a}}e^{-\frac1{|x_j|^b}}\qquad\T{ for any $a\geq0$ and for $b<1$}.
$$
Then, \eqref{3.4} (a) is obtained starting from $h^j_{z_j\bar z_j}\simgeq \frac{e^{-\frac1{|x_j|^b}}}{|x_j|^2}$, that is, the condition of type $F^2_j:=e^{-\frac1{|\delta|^b}}$ along $S_j=\R_{x_j}\setminus\{0\}$. Since $b<1$, this is infraexponential (and yields a superlogarithmic estimate for $f=\log^{\frac1b}$). \eqref{3.4} (b) follows from $|h^j_{z_j}|^2\simleq h^j_{z_j\bar z_j}$ which says that the $h^j_{z_j}$'s are not only superlogarithmic, but indeed $\frac12$-subelliptic,  multipliers.  Altogether we have that $\Box$ is hypoelliptic according to Theorem~\ref{t3.3}.

\end{document}